\documentclass[10pt]{article}

\usepackage{amssymb,amsfonts, amsopn, latexsym,graphicx,makeidx}

\newcommand{\U}{{\mathcal U}}

\newcommand{\C}{{\mathbb C}}
\newcommand{\Z}{{\mathbb Z}}

\newcommand{\D}{{\mathbb D}}

\newcommand{\mf}{{F_{f}}}
\newcommand{\mfo}{{F_{f, \mathbf 0}}}
\newcommand{\mfoo}{{F_{f_0}}}

\newcommand{\dm}{\operatorname{dim}}

\newcommand{\rank}{\operatorname{rank}}
\newcommand{\arrow}[1]{\stackrel{#1}{\longrightarrow}}
\newcommand{\text}[1]{\mbox{\rm {#1}}}

\newcommand{\pc}{{\Gamma^1_{f, z_0}}}

\newtheorem{defn0}{Definition}[section]
\newtheorem{prop0}[defn0]{Proposition}
\newtheorem{conj0}[defn0]{Conjecture}
\newtheorem{thm0}[defn0]{Theorem}
\newtheorem{lem0}[defn0]{Lemma}
\newtheorem{corollary0}[defn0]{Corollary}
\newtheorem{example0}[defn0]{Example}
\newtheorem{remark0}[defn0]{Remark}
\newtheorem{question0}[defn0]{Question}

\newenvironment{prop}{\begin{prop0}\hskip -.06in .}{\end{prop0}}

\newenvironment{thm}{\begin{thm0}\hskip -.06in .}{\end{thm0}}
\newenvironment{lem}{\begin{lem0}\hskip -.06in .}{\end{lem0}}
\newenvironment{cor}{\begin{corollary0}\hskip -.06in .}{\end{corollary0}}
\newenvironment{exm}{\begin{example0}\hskip -.06in .\rm}{\end{example0}}
\newenvironment{rem}{\begin{remark0}\hskip -.06in .\rm}{\end{remark0}}
\newenvironment{ques}{\begin{question0}\hskip -.06in .\rm}{\end{question0}}

\newcommand{\propref}[1]{Proposition~\ref{#1}}
\newcommand{\thmref}[1]{Theorem~\ref{#1}}
\newcommand{\lemref}[1]{Lemma~\ref{#1}}
\newcommand{\corref}[1]{Corollary~\ref{#1}}

\newcommand{\quesref}[1]{Question~\ref{#1}}

\newcommand{\qed}{\mbox{$\Box$}}
\newenvironment{proof}{\noindent {\bf Proof.}}{\qed\vskip 6pt}

\title{Hypersurface Singularities and the Swing\footnote{The second author would like to thank the Abdus Salam ICTP for their hospitality; most of this paper was written during a visit there.
\newline   AMS subject classifications 32B15, 32C35, 32C18, 32B10.
\newline  keywords: hypersurface singularity, Milnor fiber, swing, polar curve, vanishing cycles, discriminant, Cerf diagram, intersection diagram, perverse sheaves}}

\author{L\^e D\~ung Tr\'ang and David B. Massey}

\date{}

\begin{document}

\baselineskip= 14pt
\maketitle

\begin{abstract} Suppose that $f$ defines a singular, complex affine hypersurface. If the critical locus of $f$  is one-dimensional, we obtain  new general bounds on the ranks of the homology groups of the Milnor fiber of $f$. This result has an interesting implication on the structure of the vanishing cycles in the category of perverse sheaves.
\end{abstract}

\sloppy




\section{Introduction and Previous Results}\label{sec:intro} Let $\mathcal U$ be an open neighborhood of the origin in $\C^{n+1}$, and let $f:(\mathcal U, \mathbf 0)\rightarrow (\C, 0)$ be complex analytic. We shall always suppose that $\dm_{\mathbf 0}\Sigma f=1$, unless we explicitly state otherwise. 

\smallskip

Let $\mf =\mfo$ denote the Milnor fiber of $f$ at the origin. 
It is well-known (see \cite{katomatsu}) that the reduced integral homology, $\widetilde H_*(\mf)$, of $\mf$ can be non-zero only in degrees $n-1$ and $n$, and is free Abelian in degree $n$.  For arbitrary $f$, it is not known how to calculate, algebraically, the groups  $\widetilde H_{n-1}(\mf)$ and  $\widetilde H_n(\mf)$; in fact, it is not known how to calculate the ranks of these groups. However, there are a number of general results known for these ``top'' two homology groups of $\mf$.

\smallskip

First, we need to make some choices and establish some notation.

\smallskip

We assume that the first coordinate $z_0$ on $\U$ is a generic linear form; in the terminology of \cite{lecycles}, we need for $z_0$ to be ``prepolar'' (with respect to $f$ at the origin). This implies that, at the origin, $f_0:= f_{|_{V(z_0)}}$ has an isolated critical point, that the polar curve, $\Gamma:=\pc$, is purely $1$-dimensional at the origin (which vacuously includes the case $\Gamma=\emptyset$), and $\Gamma$ has no components contained in $V(f)$ (this last property is immediate in some definitions of the relative polar curve).

\bigskip

For convenience, we assume throughout the remainder of this paper that the neighborhood $\U$ is re-chosen, if necessary, so small that $\Sigma f\subseteq V(f)$, and every component of $\Sigma f$ and $\Gamma$ contains the origin.

\bigskip

Now, there is the attaching result of L\^e from \cite{leattach} (see, also, \cite{lecycles}), which is valid regardless of the dimension of the critical locus: 

\bigskip

\begin{thm}\label{thm:leattach} Up to diffeomorphism, $\mf$ is obtained from $\stackrel{\circ}{\D}\times\mfoo$ by attaching $\tau:=\big(\Gamma\cdot V(f))_{\mathbf 0}$ handles of index $n$.
\end{thm}

\bigskip

\begin{rem}\label{rem:leattach} On the level of homology, L\^e's attaching result is a type of Lefschetz hyperplane result; it says that, for all $i<n-1$, the inclusion map $\mfoo=\mf\cap V(z_0)\hookrightarrow\mf$ induces isomorphisms $\widetilde H_i(\mfoo)\cong \widetilde H_i(\mf)$, and $\widetilde H_n(\mf)$ and $\widetilde H_{n-1}(\mf)$ are, respectively, isomorphic to the kernel and cokernel of the boundary map $$
\Z^{\tau}\cong H_n(\mf, \mfoo)\arrow{\partial} \widetilde H_{n-1}(\mfoo)\cong\Z^{\mu_{f_0}},
$$
where $\mu_{f_0}$ denotes the Milnor number of $f_0$ at the origin. Therefore, one can certainly calculate the difference of the reduced Betti numbers of $\mf$: 
$$
\tilde b_n(\mf)-\tilde b_{n-1}(\mf) = \tau-\mu_{f_0}.
$$
Hence, bounds on one of $\tilde b_n(\mf)$ and $\tilde b_{n-1}(\mf)$ automatically produce bounds on the other.
\end{rem}

We remind the reader here of the well-known result, first proved by Teissier in \cite{teissiercargese} (in the case of an isolated singularity, but the proof works in general), that
$$
\tau \ =\  \big(\Gamma\cdot V(f))_{\mathbf 0}\  = \  \Big(\Gamma\cdot V\left(\frac{\partial f}{\partial z_0}\right)\Big)_{\mathbf 0}+\big(\Gamma\cdot V(z_0))_{\mathbf 0}.$$
As defined in \cite{lecycles}, the first summand on the right above is $\lambda^0:=\lambda^0_{f, z_0}(\mathbf 0)$, the $0$-dimensional L\^e number, and second summand on the right above is $\gamma^1:=\gamma^1_{f, z_0}(\mathbf 0)$, the $1$-dimensional polar number.

\smallskip

For each component $\nu$ of $\Sigma f$, let ${\stackrel{\circ}{\mu}}_\nu$ denote the Milnor number of $f_{|_{V(z_0-a)}}$ at a point close to the origin on $\nu\cap V(z_0-a)$, where $a$ is a small non-zero complex number. Then, 
$$
\lambda^1:=\lambda^1_{f, z_0}(\mathbf 0) :=\sum_\nu {\stackrel{\circ}{\mu}}_\nu\big(\nu\cdot V(z_0))_{\mathbf 0}
$$
is the $1$-dimensional L\^e number of $f$. Now, it is well-known, and easy to show that $\mu_{f_0}= \gamma^1+\lambda^1$. Again, see \cite{lecycles} for the above definitions and results.

\smallskip

In Proposition 3.1 of \cite{lecycles}, the second author showed how the technique of ``tilting in the Cerf diagram'' or ``the swing'', as used by L\^e and Perron in \cite{leperron} could help refine the result of \thmref{thm:leattach}. Here, we state only the homological implication of Proposition 3.1 of \cite{lecycles}.

\bigskip

\begin{thm}\label{thm:swing}  The boundary map $H_n(\mf, \mfoo)\arrow{\partial} \widetilde H_{n-1}(\mfoo)$ maps a direct summand of $H_n(\mf, \mfoo)$ of rank $\gamma^1$ isomorphically onto a direct summand of $\widetilde H_{n-1}(\mfoo)$.

Thus, the rank of $\widetilde H_n(\mf)$ is at most $\lambda^0$, and the rank of $\widetilde H_{n-1}(\mf)$ is at most $\lambda^1$.
\end{thm}

\bigskip

However, if one of the components $\nu$ of $\Sigma f$ is itself singular, then the above bounds on the ranks are known not to be optimal. A result of Siersma in \cite{siersmavarlad}, or an easy exercise using perverse sheaves (see the remark at the end of \cite{siersmavarlad}), yields:

\begin{thm}\label{thm:siersmabound} The rank of $\widetilde H_{n-1}(\mf)$ is at most $\sum_\nu {\stackrel{\circ}{\mu}}_\nu$.
\end{thm}

\bigskip

Of course, if all of the components $\nu$ of $\Sigma f$ are smooth, and $z_0$ is generic, then $\lambda^1=\sum_\nu {\stackrel{\circ}{\mu}}_\nu$, and the bounds on the ranks obtained from \thmref{thm:swing} and \thmref{thm:siersmabound} are the same. In addition,  \thmref{thm:siersmabound} is true with arbitrary field coefficients; this yields bounds on the possible torsion in $\widetilde H_{n-1}(\mf)$. We should also remark that the result of Siersma from \cite{siersmavarlad} that we cite above can actually yield a much stronger bound if one knows certain extra topological data -- specifically, one needs that the ``vertical monodromies'' are non-trivial.

\bigskip

Now, in light of \thmref{thm:swing} and \thmref{thm:siersmabound}, the question is: Is it possible that $\rank\widetilde H_{n-1}(\mf)=\lambda^1$?

\bigskip

Of course, the answer to this question is ``yes''; if $f$ has a smooth critical locus which defines a family of isolated singularities with constant Milnor number $\mu_{f_0}$, then certainly $\widetilde H_{n}(\mf)=0$ and $\widetilde H_{n-1}(\mf)\cong\Z^{\lambda^1}=\Z^{\mu_{f_0}}$. We refer to this case as the {\bf trivial case}. It is important to note that being in the trivial case implies that $V(z_0)$ transversely intersects the smooth critical locus at the origin.

\bigskip

By the non-splitting result, proved independently by Gabrielov \cite{gabrielov}, Lazzeri \cite{lazzeri}, and L\^e \cite{leacampo}, we have:

\smallskip

\begin{prop}\label{prop:trivequiv} The trivial case is equivalent to the case $\Gamma=\emptyset$.
\end{prop}

\bigskip

Now, we can state our Main Theorem:

\bigskip

\noindent{\bf Main Theorem}. {\it Suppose that $\dm_{\mathbf 0}\Sigma f=1$ and $\dm_{\mathbf 0}\Sigma f_0=0$. Then, the following are equivalent:

\vskip .1in

\noindent a) We are in the trivial case, i.e., $f$ has a smooth critical locus which defines a family of isolated singularities with constant Milnor number $\mu_{f_0}$;

\vskip .1in

\noindent b) $\rank\widetilde H_{n-1}(\mf)=\lambda^1$;

\vskip .1in

\noindent c) there exists a prime $p$ such that $\dm \widetilde H_{n-1}(\mf;\ \Z/p\Z)=\lambda^1$.

\vskip .1in

Thus, if we are not in the trivial case, $\rank\widetilde H_{n-1}(\mf)<\lambda^1$, and so $\rank\widetilde H_{n}(\mf)<\lambda^0$, and these inequalities hold with $\Z/p\Z$ coefficients (here, $p$ is prime).}

\bigskip

\begin{rem}\label{rem:mainthm} We remark again that if one of the components of $\Sigma f$ is itself singular (and, hence, we are not in the trivial case), then the conclusion that $\rank\widetilde H_{n-1}(\mf)<\lambda^1$ already follows from \thmref{thm:siersmabound}. Even in the case where all of the components of $\Sigma f$ are smooth, we could conclude that $\rank\widetilde H_{n-1}(\mf)<\lambda^1$ from \cite{siersmavarlad} {\bf if} we knew that one of the vertical monodromies were non-trivial. However, the vertical monodromies are fairly complicated topological data to calculate, and it is also true that the vertical monodromies can be trivial even when the polar curve is non-empty, i.e., when we are not in the trivial case. Thus, our Main Theorem cannot be proved by analyzing the vertical monodromies.

In \cite{siersmaisoline}, Siersma proved another closely related result. On the level of homology, what he proved was that, if we are not in the trivial case, and $\Sigma f$ has a single smooth component, $\nu$, such that ${\stackrel{\circ}{\mu}}_\nu =1$, then $\widetilde H_{n-1}(\mf)=0$; our Main Theorem, including the modulo $p$ statement, is a strict generalization of this.

In addition, we should point out that, in \cite{dejong}, Th. de Jong provides evidence that a result like our Main Theorem might be true. 
\end{rem}

\bigskip

We prove our Main Theorem by combining the swing technique of \thmref{thm:swing} and the connectivity of the vanishing cycle intersection diagram for isolated singularities, as was proved independently by Gabrielov in \cite{gabrielov} and Lazzeri in \cite{lazzeri}. In some recent notes, M. Tib\u ar uses similar techniques and reaches a number of conclusions closely related to our result.

\bigskip

As a corollary to our Main Theorem, we show that it implies that the vanishing cycles of $f$, as an object in the category of perverse sheaves, cannot be semi-simple in non-trivial cases where $\Sigma f$ has smooth components of arbitrary dimension.

\bigskip

In the final section of this paper, we make some final remarks and present counterexamples to some conceivable ``improvements'' on the statement of the Main Theorem.

\section{The Swing}\label{sec:swing}

In the Introduction, we referred to the swing (or, tilting in the Cerf diagram), which was used by L\^e and Perron in \cite{leperron} and in Proposition 3.1 of \cite{lecycles}, where the swing was used to prove \thmref{thm:swing}. The swing has also been studied in \cite{caubelthesis}, \cite{tibarbouq}, \cite{lecycles}, \cite{vannierthesis}. As the swing is so crucial to the proof of the main theorem, we wish to give a careful explanation of its construction.

\bigskip

Suppose that $\mathcal W$ is an open neighborhood of the origin in $\mathbb C^2$. We will use the coordinates $x$ and $y$ on $\mathcal W$. For notational ease, when we restrict $x$ and $y$ to various subspaces where the domain is clear, we shall continue to write simply $x$ and $y$.

 Let $C$ be a complex analytic curve in $\mathcal W$ such that every component of $C$ contains the origin. We assume that the origin is an isolated point in $V(x)\cap C$ and in $V(y)\cap C$, i.e., that $C$ does not have a component along the $x$- or $y$-axis.
 
 Below, we let $\mathbb D_\epsilon$ denote a closed disk, of radius $\epsilon$, centered at the origin, in the complex plane. We denote the interior of $\mathbb D_\epsilon$ by ${\stackrel{\circ}{\mathbb D}}_{\epsilon}$, and when we delete the origin, we shall superscript with an asterisk, i.e., $\mathbb D^*_\epsilon:= \mathbb D_\epsilon-\{0\}$ and 
 ${\stackrel{{}\hskip -.06in \circ}{\mathbb D^*_{\epsilon}}}:={\stackrel{\circ}{\mathbb D}}_{\epsilon}-\{0\}$.

\bigskip

We select $0<\epsilon_2\ll\epsilon_1\ll 1$ so that:

\medskip

\noindent i): the ``half-open'' polydisk $\mathbb D_{\epsilon_1}\times{\stackrel{\circ}{\mathbb D}}_{\epsilon_2}$ is contained in $\mathcal W$;

 \medskip
 
 \noindent ii):   $(\partial \mathbb D_{\epsilon_1}\times{\stackrel{\circ}{\mathbb D}}_{\epsilon_2})\cap C=\emptyset$ (this uses that the origin is an isolated point in $V(y)\cap C$) ;
 
 \medskip
 
 Note that ii) implies that $(\mathbb D_{\epsilon_1}\times{\stackrel{\circ}{\mathbb D}}_{\epsilon_2})\cap C = ({\stackrel{\circ}{\mathbb D}}_{\epsilon_1}\times{\stackrel{\circ}{\mathbb D}}_{\epsilon_2})\cap C$.
 
 \medskip
 
 \noindent iii):  $\mathbb D_{\epsilon_1}\times {\stackrel{{}\hskip -.06in \circ}{\mathbb D^*_{\epsilon_2}}}\arrow{y} {\stackrel{{}\hskip -.06in \circ}{\mathbb D^*_{\epsilon_2}}}$ is a proper stratified submersion, where the Whitney strata are $\partial \mathbb D_{\epsilon_1}\times {\stackrel{{}\hskip -.06in \circ}{\mathbb D^*_{\epsilon_2}}}$, \hbox{$({\stackrel{\circ}{\mathbb D}}_{\epsilon_1}\times {\stackrel{{}\hskip -.06in \circ}{\mathbb D^*_{\epsilon_2}}})-C$}, and $({\stackrel{\circ}{\mathbb D}}_{\epsilon_1}\times {\stackrel{{}\hskip -.06in \circ}{\mathbb D^*_{\epsilon_2}}})\cap C$.
 
 \medskip
 
 \noindent iv): $({\stackrel{\circ}{\mathbb D}}_{\epsilon_1}\times {\stackrel{{}\hskip -.06in \circ}{\mathbb D^*_{\epsilon_2}}})\cap C\arrow{y} {\stackrel{{}\hskip -.06in \circ}{\mathbb D^*_{\epsilon_2}}}$ is an $m$-fold covering map, where $m:=(C\cdot V(y))_\mathbf 0$.

\bigskip

Let $D:= (\mathbb D_{\epsilon_1}\times{\stackrel{\circ}{\mathbb D}}_{\epsilon_2})\cap (C\cup V(y))$. Let $(x_0, y_0)\in (\mathbb D^*_{\epsilon_1}\times {\stackrel{{}\hskip -.06in \circ}{\mathbb D^*_{\epsilon_2}}})-D$. Let $\sigma:[0,1]\rightarrow\{x_0\}\times {\stackrel{\circ}{\mathbb D}}_{\epsilon_2}$ be a smooth, simple path such that $\sigma(0)= (x_0, y_0)$, $\sigma(1)=:(x_0, y_1)\in C$, and $\sigma([0,1))\subseteq (\{x_0\}\times{\stackrel{\circ}{\mathbb D}}_{\epsilon_2})-D$. 

Let $S$ be the image of $\sigma$; as $\sigma$ is simple, $S$ is homeomorphic to $[0, 1]$. Let $\sigma_0:=y\circ\sigma$ and let $S_0$ be the image of $\sigma_0$. Thus, $S_0$ is homeomorphic to $[0,1]$ and is contained in ${\stackrel{{}\hskip -.06in \circ}{\mathbb D^*_{\epsilon_2}}}$.

\bigskip

\begin{lem}\label{lem:swing}{\rm({\bf The Swing})}  There exists a continuous function $H:[0,1]\times [0,1]\rightarrow \mathbb D_{\epsilon_1}\times S_0$ with the following properties:

\smallskip

\noindent a) $H(t, 0)=\sigma(t)$, for all $t\in[0,1]$;

\smallskip

\noindent b) $H(t, 1)\in \mathbb D_{\epsilon_1}\times\{y_0\}$, for all $t\in[0,1]$;

\smallskip

\noindent c) $H(0, u) = (x_0, y_0)$

\smallskip

\noindent d) if $H(t,u)\in D$, then $t=1$;

\smallskip

\noindent e) $H(1,u)\in C$, for all $u\in[0,1]$;

\smallskip

\noindent f) the path $\eta$ given by $\eta(u):=H(1,u)$ is a homeomorphism onto its image.
\medskip

Thus, $H$ is a homotopy from $\sigma$ to the path $\gamma$ given by $\gamma(t):=H(t, 1)\in \mathbb D_{\epsilon_1}\times\{y_0\}$, such that $(x_0, y_0)$ is ``fixed'' and the point $(x_0, y_1)=H(1, 0)$ ``swings up to the point'' $H(1,1)$ by ``sliding along'' $C$, while the remainder of $\sigma$ does not hit $D$ as it ``swings up'' to $\gamma$.

\end{lem}
\begin{proof} The proper stratified submersion $\mathbb D_{\epsilon_1}\times {\stackrel{{}\hskip -.06in \circ}{\mathbb D^*_{\epsilon_2}}}\arrow{y} {\stackrel{{}\hskip -.06in \circ}{\mathbb D^*_{\epsilon_2}}}$ is a locally trivial fibration, where the local trivialization respects the strata. The restriction of this fibration $\mathbb D_{\epsilon_1}\times S_0\arrow{y} S_0$ is a locally trivial fibration over a contractible space and, hence, is equivalent to the trivial fibration. 

Therefore, there exists a homeomorphism 
$$\Psi: \big(\mathbb D_{\epsilon_1}\times S_0, (\mathbb D_{\epsilon_1}\times S_0)\cap C\big)\rightarrow \big(\mathbb D_{\epsilon_1}\times \{y_0\}, (\mathbb D_{\epsilon_1}\times \{y_0\})\cap C\big)\times[0,1]$$
such that the projection of $\Psi(x, \sigma_0(t))$ onto the $[0,1]$ factor is simply $t$, and such that $\Psi(x, y_0)= ((x, y_0), 0)$. It follows that there is a path $\alpha:[0,1]\rightarrow \mathbb D_{\epsilon_1}$ such that $\Psi(\sigma(t)) = ((\alpha(t), y_0), t)$, for all $t\in[0,1]$. Define $H:[0,1]\times [0,1]\rightarrow \mathbb D_{\epsilon_1}\times S_0$ by
$$
H(t, u):= \Psi^{-1}\big((\alpha(t), y_0), (1-u)t\big).
$$
All of the given properties of $H$ are now trivial to verify.
\end{proof}

\bigskip

\begin{rem}\label{rem:triangle} By Property c) of \lemref{lem:swing}, the map $H$ yields a corresponding map $H^T$ whose domain is a triangle instead of a square. One pictures the image of $H$, or of $H^T$, as a ``gluing in'' of this triangle into $\mathbb D_{\epsilon_1}\times S_0$ in such a way that one edge of the triangle is glued diffeomorphically to $S$, and another edge is glued diffeomorphically onto the image of $\eta$. The third edge of the triangle is glued onto the image of $\gamma$, but not necessarily in a one-to-one fashion.
\end{rem}

\section{The Main Theorem}\label{sec:main} 

In this section, we will prove the Main Theorem, as stated in the Introduction and as appears below as \thmref{thm:main}. That a) of the Main Theorem implies both b) and c) is well-known; one can, for instance, conclude it from \thmref{thm:leattach}. The difficulty is to prove that b) and c) imply a). In fact, we prove the contrapositives; we prove that if we are not in the trivial case, then $\rank\widetilde H_{n-1}(\mf)<\lambda^1$ and $\dm\widetilde H_{n-1}(\mf;\ \mathbb Z/p\mathbb Z)<\lambda^1$. We must first describe the machinery that goes into this part of the proof.

\bigskip

As the value of $\lambda^1$ is minimal for generic $z_0$, we lose no generality if we assume that our linear form $z_0$ is chosen more generically than simply being prepolar. We choose $z_0$ so generically that, in addition to being prepolar, the discriminant, $D$, of the map $(z_0, f)$ and the corresponding Cerf diagram, $C$, have the usual properties -- as given, for instance, in \cite{leperron}, \cite{tibarbouq}, and \cite{vannierthesis}. We will describe the needed properties below.

\bigskip

Let $\widetilde \Psi:=(z_0, f):(\U, \mathbf 0)\rightarrow (\C^2, \mathbf 0)$. We use the coordinates $(u, v)$ on $\C^2$. The critical locus $\Sigma \widetilde \Psi$ of $\widetilde \Psi$ is the union of $\Sigma f$ and $\Gamma$. The discriminant $D:=\widetilde \Psi(\Sigma \widetilde \Psi)$ consists of the $u$-axis together with the Cerf diagram $C:=\overline{D-V(v)}$. We assume that $z_0$ is generic enough so that the polar curve is reduced and that, in a neighborhood of the origin, $\widetilde \Psi_{|_\Gamma}$ is one-to-one.

We choose real numbers $\epsilon$, $\delta$, and $\omega$ so that $0< \omega\ll \delta\ll\epsilon\ll 1$. Let $B_\epsilon\subseteq\C^n$ be a closed ball, centered at the origin, of radius $\epsilon$. Let ${\stackrel{\circ}{\D}}_\delta$ and ${\stackrel{\circ}{\D}}_\omega$ be open disks in $\C$, centered at $0$, of radii $\delta$ and $\omega$, respectively.

One considers the map from $({\stackrel{\circ}{\D}}_\delta\times B_\epsilon)\cap f^{-1}({\stackrel{\circ}{\D}}_\omega)$ onto ${\stackrel{\circ}{\D}}_\delta\times{\stackrel{\circ}{\D}}_\omega$ given by the restriction of $\widetilde  \Psi$; we let $ \Psi$ denote this restriction. As $B_\epsilon$ is a closed ball, the map $ \Psi$ is certainly proper, but the domain has an interior stratum, and a stratum coming from the boundary of $B_\epsilon$. However, for generic $z_0$, all of the stratified critical points lie on $\Gamma\cup\Sigma f$, i.e., above $D$.  

We continue to write simply $D$ and $C$, in place of $D\cap({\stackrel{\circ}{\D}}_\delta\times{\stackrel{\circ}{\D}}_\omega)$ and $C\cap({\stackrel{\circ}{\D}}_\delta\times{\stackrel{\circ}{\D}}_\omega)$. As $ \Psi$ is a proper stratified submersion above ${\stackrel{\circ}{\D}}_\delta\times{\stackrel{\circ}{\D}}_\omega- D$, and as $ \Psi_{|_\Gamma}$ is one-to-one, many homotopy arguments in $({\stackrel{\circ}{\D}}_\delta\times B_\epsilon)\cap f^{-1}({\stackrel{\circ}{\D}}_\omega)$ can be obtained from lifting constructions in ${\stackrel{\circ}{\D}}_\delta\times{\stackrel{\circ}{\D}}_\omega$. This is the point of considering the discriminant and Cerf diagram. 

Let $v_0\in{\stackrel{\circ}{\D}}_\omega-\{\mathbf 0\}$. By construction, up to diffeomorphism, $ \Psi^{-1}({\stackrel{\circ}{\D}}_\delta\times\{v_0\})$ is $\mf$ and $ \Psi^{-1}((0, v_0))$ is $\mfoo$. In fact, for all $u_0$, where $|u_0|\ll |v_0|$, $ \Psi^{-1}((u_0, v_0))$ is diffeomorphic to $\mfoo$; we fix such a non-zero $u_0$, and let $\mathbf a:=(u_0, v_0)$.

\smallskip

We wish to pick a distinguished basis for the vanishing cycles of $f_0$ at the origin, as in I.1 of \cite{agv} (see, also, \cite{dimcasing}). We do this by selecting paths in $\{u_0\}\times {\stackrel{\circ}{\D}}_\omega$ which originate at $\mathbf a$. We must be slightly careful in how we do this. 

First, fix a path $p_0$ from $\mathbf a$ to $(u_0, 0)$. Select paths  $q_1, \dots, q_{\gamma^1}$ from $\mathbf a$ to each of the points in $(\{u_0\}\times {\stackrel{\circ}{\D}}_\omega)\cap C=:\{y_1, \dots, y_{\gamma^1}\}$. The paths $p_0, q_1, \dots, q_{\gamma^1}$ should not intersect each other and should intersect the set $\{(u_0, 0), y_1, \dots, y_{\gamma^1}\}$  only at the endpoints of the paths. Moreover, when at the point $\mathbf a$, the paths $p_0, q_1, \dots, q_{\gamma^1}$ should be in clockwise order. Let $r_0$ be a clockwise loop very close to $p_{0}$, from $\mathbf a$ around $(u_0, 0)$.

As we are not assuming that $f$ had an isolated line singularity, we must perturb $f_{|_{V(z_0-u_0)}}$ slightly to have $(u_0, 0)$ split into $\lambda^1$ points, $x_1, \dots, x_{\lambda^1}$ inside the loop $r_0$; each of these points corresponds to an $A_1$ singularity in the domain. We select paths $p_1, \dots, p_{\lambda^1}$ from $\mathbf a$ to each of the points $x_1, \dots, x_{\lambda^1}$, and paths $q_1, \dots, q_{\gamma^1}$ from $\mathbf a$ to each of the points in $(\{u_0\}\times {\stackrel{\circ}{\D}}_\omega)\cap C=:\{y_1, \dots, y_{\gamma^1}\}$. We may do this in such a way that the paths $p_1, \dots, p_{\lambda^1}, q_1, \dots, q_{\gamma^1}$ are in clockwise order. 

\smallskip

The lifts of these paths via the perturbed $f_{|_{V(z_0 - u_0)}}$ yield representatives of elements of $H_{n+1}(B_\epsilon, \mfoo)$, whose boundaries in $\widetilde H_n(\mfoo)$ form a distinguished basis  $\Delta_1^\prime, \dots, \Delta_{\lambda^1}^\prime, \Delta_1, \dots, \Delta_{\gamma^1}$.

\smallskip

By using the swing (\lemref{lem:swing}), the paths $q_1, \dots, q_{\gamma^1}$ are taken to new paths $\hat q_1, \dots, \hat q_{\gamma^1}$ in ${\stackrel{\circ}{\D}}_\delta\times\{v_0\}$. Each $\hat q_i$ path represents a relative homology class in $H_n(\mf, \mfoo)$ whose boundary in $\widetilde H_{n-1}(\mfoo)$ is precisely $\Delta_i$.  \thmref{thm:swing} follows from this. 

\bigskip

We can now prove the Main Theorem:

\bigskip

\begin{thm}\label{thm:main} Suppose that $\dm_{\mathbf 0}\Sigma f=1$ and $\dm_{\mathbf 0}\Sigma f_0=0$. Then, the following are equivalent:

\vskip .1in

\noindent a) We are in the trivial case, i.e., $f$ has a smooth critical locus which defines a family of isolated singularities with constant Milnor number $\mu_{f_0}$;

\vskip .1in

\noindent b) $\rank\widetilde H_{n-1}(\mf)=\lambda^1$;

\vskip .1in

\noindent c) there exists a prime $p$ such that $\dm \widetilde H_{n-1}(\mf;\ \Z/p\Z)=\lambda^1$.

\vskip .1in

Thus, if we are not in the trivial case, $\rank\widetilde H_{n-1}(\mf)<\lambda^1$, and so $\rank\widetilde H_{n}(\mf)<\lambda^0$, and these inequalities hold with $\Z/p\Z$ coefficients (here, $p$ is prime).\end{thm}

\begin{proof} As mentioned above, that a) implies b) and c) is well-known. Assume then that we are not in the trivial case. We will prove that $\rank\widetilde H_{n-1}(\mf)<\lambda^1$, and then indicate why the same proof applies with $\mathbb Z/p\mathbb Z$ coefficients.

By \propref{prop:trivequiv}, $\Gamma\neq\emptyset$, and so $C\neq\emptyset$. We want to construct just one new path in  $\{u_0\}\times {\stackrel{\circ}{\D}}_\omega$, one which originates at $\mathbf a$, ends at a point of $C$, and misses all of the other points of $D$; we want this path to swing up to a path in  ${\stackrel{\circ}{\D}}_\delta\times\{v_0\}$, and represent a relative homology class in $H_n(\mf, \mfoo)$ whose boundary is not in the span of $\Delta_1, \dots, \Delta_{\gamma^1}$.

By the connectivity of the vanishing cycle intersection diagram (\cite{gabrielov}, \cite{lazzeri}), one of the $\Delta^\prime_j$ must have a non-zero intersection pairing with one of the $\Delta_i$, i.e., there exist $i_0$ and $j_0$ such that $\langle \Delta_{i_0}, \Delta^\prime_{j_0}\rangle\neq 0$. 

By fixing the path $p_{j_0}$ and all the $q_i$ paths, but reselecting the other $p_j$, for $j\neq j_0$, we may assume that $j_0=1$, i.e., that $\langle \Delta_{i_0}, \Delta^\prime_{1}\rangle\neq 0$. 

We follow now Chapter 3.3 of \cite{dimcasing}. Associated to each path $p_{j}$, $1\leqslant j\leqslant\lambda^1$, is a (partial) monodromy automorphism $T^\prime_{j}:\widetilde H_{n-1}(\mfoo)\rightarrow \widetilde H_{n-1}(\mfoo)$, induced by taking a clockwise loop $r_{j}$ very close to $p_{j}$, from $\mathbf a$ around $x_{j}$. Let $T^\prime:= T^\prime_1\dots T^\prime_{\lambda^1}$, where composition is written in the order of \cite{dimcasing}. We claim that $T^\prime(\Delta_{i_0})$ is in the image of $\delta: H_n(\mf, \mfoo)\rightarrow H_{n-1}(\mfoo)$, but is not in $\operatorname{Span}\{\Delta_1, \dots, \Delta_{\gamma^1}\}$.

The composition $r$ of the loops $r_1, \dots, r_{\lambda^1}$ is homotopy-equivalent, in $\{u_0\}\times {\stackrel{\circ}{\D}}_\omega-\big\{\{x_1, \dots, x_{\lambda^1}\}\cup C\big\}$, to the loop $r_0$ (from our discussion before the theorem).   By combining (concatenating) the loop $r_{0}$ and the path $q_{i_0}$, we obtain a path in $\{u_0\}\times {\stackrel{\circ}{\D}}_\omega$ which is homotopy-equivalent to a simple path which swings up to a corresponding path in  ${\stackrel{\circ}{\D}}_\delta\times\{v_0\}$. Thus, $T^\prime(\Delta_{i_0})$ is in the image of $\delta$.

Now, by the Corollaries to the Picard-Lefschetz Theorem in \cite{agv}, p. 26, or as in \cite{dimcasing}, Formula 3.11,
$$
T^\prime(\Delta_{i_0}) = \Delta_{i_0}- (-1)^{\frac{n(n-1)}{2}}\langle \Delta_{i_0}, \Delta^\prime_{1} \rangle \Delta^\prime_{1}+\beta_2 \Delta^\prime_{2}+\dots+\beta_{\lambda^1} \Delta^\prime_{\lambda^1},
$$
for some integers $\beta_2, \dots, \beta_{\lambda^1}$.
As the $\Delta_1^\prime, \dots, \Delta_{\lambda^1}^\prime, \Delta_1, \dots, \Delta_{\gamma^1}$ form a basis, and as $\langle \Delta_{i_0}, \Delta^\prime_{1} \rangle\neq 0$, $T^\prime(\Delta_{i_0})$ is not in $\operatorname{Span}\{\Delta_1, \dots, \Delta_{\gamma^1}\}$.

This finishes the proof over the integers. Over $\Z/p\Z$, the proof is identical, since the intersection diagram is also connected modulo $p$; see \cite{gabrielov}.
\end{proof}

\bigskip

\begin{rem}\label{rem:careful} One must be careful in the proof above; it is tempting to try to use simply $T^\prime_1(\Delta_{i_0})$ in place of $T^\prime(\Delta_{i_0})$. The problem with this is that $T^\prime_1(\Delta_{i_0})$ is not represented by a path in $\{u_0\}\times {\stackrel{\circ}{\D}}_\omega - \{(u_0, 0)\}$ and, thus, there is no guaranteed swing isotopy to a corresponding path in  ${\stackrel{\circ}{\D}}_\delta\times\{v_0\}$.
\end{rem}

\bigskip

In the corollary below, we obtain a conclusion when the dimension of $\Sigma f$ is arbitrary. We use the notation and terminology from \cite{lecycles}. In particular, $\lambda^s_{f, \mathbf z}(\mathbf 0)$ is the $s$-dimensional L\^e number of $f$ at the origin with respect to the coordinates $\mathbf z$.

\bigskip

\begin{cor}\label{cor:lebound} Suppose that the dimension of $\Sigma f$ at the origin is $s$, where $s\geqslant 1$ is arbitrary. Assume that the coordinates $\mathbf z:=(z_0, ..., z_{s-1})$ are prepolar for $f$ at the origin, and that the $s$-dimensional relative polar variety $\Gamma^s_{f, \mathbf z}$ at the origin is not empty.

Then, both $\rank\widetilde H_{n-s}(\mf)$ and
$\dm \widetilde H_{n-s}(\mf;\ \Z/p\Z)$ are strictly less than $\lambda^s_{f, \mathbf z}(\mathbf 0)$.
\end{cor}
\begin{proof} One simply takes the codimension $s-1$ linear slice $N:=V(z_0, \dots, z_{s-2})$ through the origin. Then, $f_{|_N}$ has a $1$-dimensional critical locus and, by iterating \thmref{thm:leattach}, $\widetilde H_{n-s}(\mf)\cong \widetilde H_{(n-s+1)-1}(F_{f_{|_N}})$. Now, by Proposition 1.21 of \cite{lecycles}, $\lambda^s_{f, \mathbf z}(\mathbf 0) = \lambda^1_{f_{|_N}, z_{s-1}}(\mathbf 0)$. The corollary nows follows at once from \thmref{thm:main} (the proof with $\Z/p\Z$ coefficients is identical).
\end{proof}

\bigskip

As we shall see, \corref{cor:lebound} puts restrictions on the types of perverse sheaves that one may obtain as vanishing cycles of the shifted constant sheaf on affine space. Below, we refer to the constant sheaf on $\nu$ of dimension ${\stackrel{\circ}{\mu}}_\nu$, shifted by $1$ and extended by zero to all of $V(f)$; we write $(k^{{\stackrel{\circ}{\mu}}_\nu})^\bullet_\nu[1]$ for this sheaf (note that we omit the reference to the extension by zero in the notation). The isomorphisms and direct sums that we write below are in the Abelian category of perverse sheaves.

\bigskip

In the trivial case, $\Sigma f$ consists of a single smooth component $\nu$ and $\phi_f[-1]k^\bullet_{\U}[n+1]\cong (k^{{\stackrel{\circ}{\mu}}_\nu})^\bullet_\nu[1]$. Aside from the trivial case, is it possible for $(k^{{\stackrel{\circ}{\mu}}_\nu})^\bullet_\nu[1]$ to be a direct summand of $\phi_f[-1]k^\bullet_{\U}[n+1]$? The following corollary provides a partial answer, and generalizes the question/answer to critical loci of arbitrary dimension.

\bigskip

\begin{cor}\label{cor:pervstruct} Suppose that the critical locus of $f$ is  $s$-dimensional, where $s\geqslant 1$ is arbitrary. For each $s$-dimensional component $\nu$ of $\Sigma f$, let ${\stackrel{\circ}{\mu}}_\nu$ denote the Milnor number of $f$ restricted to a generic normal slice of $\nu$.

If  $\Sigma f$ is smooth  and the generic $s$-dimensional relative polar variety of $f$ is empty, then  $\phi_f[-1]k^\bullet_{\U}[n+1]\cong (k^{{\stackrel{\circ}{\mu}}_\nu})^\bullet_\nu[s]$.

If each component of $\Sigma f$ is smooth, and the generic $s$-dimensional relative polar variety of $f$ is not empty, then  $\bigoplus_\nu(k^{{\stackrel{\circ}{\mu}}_\nu})^\bullet_\nu[s]$ is not a direct summand of $\phi_f[-1]k^\bullet_{\U}[n+1]$.
\end{cor}

\begin{proof}
If $\Sigma f$ is smooth and the $s$-dimensional relative polar variety is empty, $V(f)$ has an $a_f$ stratification consisting of two strata: $V(f)-\Sigma f$ and $\Sigma f$. As $\phi_f[-1]k^\bullet_{\U}[n+1]$ is constructible with respect to any $a_f$ stratification, the first statement follows.

\smallskip

If each component of $\Sigma f$ is smooth, then, for generic coordinates, the $s$-dimensional L\^e number $\lambda^s_f(\mathbf 0)$ will be equal to $\sum_\nu {\stackrel{\circ}{\mu}}_\nu$, where we sum over $s$-dimensional components. Now, the second statement follows at once from \corref{cor:lebound}, since such a direct summand would immediately imply that the dimension of $\widetilde H_{n-s}(\mf)$ is too big.
\end{proof}

\section{Comments, Questions, and Counterexamples}\label{sec:cqc}

 One might hope that a stronger result than \thmref{thm:main} is true.
 
\bigskip

For instance, given that \thmref{thm:main} and \thmref{thm:siersmabound} are true, it is natural to ask the following:

\begin{ques} If we are not in the trivial case, is the rank of $\widetilde H_{n-1}(\mf)$ strictly less than $\sum_\nu {\stackrel{\circ}{\mu}}_\nu$?
\end{ques}

\smallskip

The answer to the above question is ``no''. One can find examples of this in the literature, but perhaps the easiest is the following:

\smallskip

\begin{exm} Let $f:= (y^2-x^3)^2+w^2$. Then, $\Sigma f$ has a single component $\nu:=V(w, y^2-x^3)$, and one easily checks that ${\stackrel{\circ}{\mu}}_\nu = 1$. However, as $f$ is the suspension of $(y^2-x^3)^2$, the Sebastiani-Thom Theorem (here, we need the version proved by Oka in \cite{okasebthom}) implies
$$
\widetilde H_1(\mf)\cong\widetilde H_0(F_{(y^2-x^3)^2})\cong\Z.
$$
Moreover, by suspending $f$ again, one may produce an example in which $f$ itself has a single irreducible component at the origin.
\end{exm}

\bigskip

 Now, let $\alpha$ be the number of irreducible components of $\Sigma f$.
 
 \smallskip
 
 \begin{ques}\label{ques:betterbound} If we are not in the trivial case, is the rank of $\widetilde H_{n-1}(\mf)$ strictly less than $\lambda^1-\alpha$?
\end{ques}

\smallskip

Again, there are many examples in the literature which demonstrate that the answer to this question is ``no''. One simple example is:

\begin{exm} 
The function $f=x^2y^2+w^2$ has a critical locus consisting of two lines, $\lambda^1=2$, but -- using the Sebastiani-Thom Theorem again -- we find that $\widetilde H_1(\mf)\cong\Z$.
\end{exm}

\bigskip

However, a result such as that asked about in \quesref{ques:betterbound}, but where $\alpha$ is replaced by a quantity involving the number of components of $\Gamma$, or numbers of various types of components in the Cerf diagram, seems more likely. Moreover, if we put more conditions on the intersection diagram for the vanishing cycles of $f_0$, we could certainly obtain sharper bounds than we do in the Main Theorem.  Or, if we know more topological data, such as the vertical monodromies, as in \cite{siersmavarlad}, we could obtain better bounds. However, other than \thmref{thm:main} , we know of no nice, effectively calculable, formula which holds in all cases.

\bigskip

Finally, \corref{cor:pervstruct} leads us to ask: 

\smallskip

\begin{ques}\label{ques:pervstruct} Which perverse sheaves can be obtained as the vanishing cycles of the constant sheaf on affine space?
\end{ques}

\bigskip

Unlike our previous questions, we do not know the answer to \quesref{ques:pervstruct}.

\newpage
\bibliographystyle{plain}
\bibliography{Masseybib}

\begin{thebibliography}{10}

\bibitem{agv}
{Arnold, V. I., Gusein-Zade, S. M., Varchenko, A. N.}
\newblock {\em {Singularities of Differentiable Maps II, Monodromy and
  Asymptotics of Integrals}}.
\newblock {Birkh\"auser}, 1988.

\bibitem{caubelthesis}
{Caubel, C.}
\newblock {\em {Sur la topologie d'une famille de pinceaux de germes
  d'hypersurfaces complexes}}.
\newblock PhD thesis, Universit\'e Toulouse III, 1998.

\bibitem{dejong}
{de Jong, Th.}
\newblock {Some classes of line singularities}.
\newblock {\em Math. Zeitschrift}, 198:493--517, 1998.

\bibitem{dimcasing}
{Dimca, A.}
\newblock {\em {Singularities and Topology of Hypersurfaces}}.
\newblock Universitext. Springer-Verlag, 1992.

\bibitem{gabrielov}
{Gabrielov, A. M.}
\newblock {Bifurcations, Dynkin Diagrams, and Modality of Isolated
  Singularities}.
\newblock {\em Funk. Anal. Pril.}, 8 (2):7--12, 1974.

\bibitem{katomatsu}
{Kato, M. and Matsumoto, Y.}
\newblock {On the connectivity of the Milnor fibre of a holomorphic function at
  a critical point}.
\newblock {\em Proc. of 1973 Tokyo manifolds conf.}, pages 131--136, 1973.

\bibitem{lazzeri}
{Lazzeri, F.}
\newblock {\em {Some Remarks on the Picard-Lefschetz Monodromy}}.
\newblock {Quelques journ\'ees singuli\`eres}. Centre de Math. de l'Ecole
  Polytechnique, Paris, 1974.

\bibitem{leacampo}
{L\^e, D. T. }.
\newblock {Une application d'un th\'eor\`eme d'A'Campo a l'equisingularit\'e}.
\newblock {\em Indag. Math}, 35:403--409, 1973.

\bibitem{leattach}
{L\^e, D. T.}
\newblock {Calcul du Nombre de Cycles Evanouissants d'une Hypersurface
  Complexe}.
\newblock {\em Ann. Inst. Fourier, Grenoble}, 23:261--270, 1973.

\bibitem{leperron}
{L\^e, D. T. and Perron, B.}
\newblock {Sur la fibre de Milnor d'une singularit\'e isol\'ee en dimension
  complexe trois}.
\newblock {\em C. R. Acad. Sci. Pairs S\'er. A}, 289:115--118, 1979.

\bibitem{lecycles}
{Massey, D.}
\newblock {\em {L\^e Cycles and Hypersurface Singularities}}, volume 1615 of
  {\em Lecture Notes in Math.}
\newblock Springer-Verlag, 1995.

\bibitem{okasebthom}
{Oka, M.}
\newblock {On the homotopy type of hypersurfaces defined by weighted
  homogeneous polynomials}.
\newblock {\em Topology}, 12:19--32, 1973.

\bibitem{siersmaisoline}
{Siersma, D.}
\newblock {Isolated line singularities}.
\newblock {\em Proc. Symp. Pure Math.}, 35, part 2, Arcata Singularities
  Conf.:485--496, 1983.

\bibitem{siersmavarlad}
{Siersma, D.}
\newblock {Variation mappings on singularities with a $1$-dimensional critical
  locus}.
\newblock {\em Topology}, 30:445--469, 1991.

\bibitem{teissiercargese}
{Teissier, B.}
\newblock {Cycles \'evanescents, sections planes et conditions de Whitney}.
\newblock {\em Ast\'erisque}, 7-8:285--362, 1973.

\bibitem{tibarbouq}
{Tib\u ar, M.}
\newblock {Bouquet Decomposition of the Milnor Fiber}.
\newblock {\em Topology}, 35:227--241, 1996.

\bibitem{vannierthesis}
{Vannier, J. P.}
\newblock {\em {Familles \`a param\`etre de fonctions holomorphes \`a ensemble
  singulier de dimension z\'ero ou un}}.
\newblock PhD thesis, Dijon, 1987.

\end{thebibliography}
\end{document}